\documentclass[11pt, a4paper]{amsart}
\usepackage{amssymb, array, amsmath, amscd, pdfpages, enumerate, amsthm, fixltx2e, setspace, hyperref,mathtools}

\usepackage[utf8]{inputenc}
\usepackage[english]{babel}
\usepackage{microtype}
\usepackage{ifthen}
\usepackage{color}
\usepackage{enumerate}

\usepackage{graphics}
\usepackage{epsfig}
\usepackage{amssymb} 
\usepackage{amsmath}
{
		\newtheorem{theorem}{Theorem}
       \newtheorem{assumption}[theorem]{Assumption}
       \newtheorem{definition}[theorem]{Definition}
       \newtheorem{lemma}[theorem]{Lemma}
	   
}

\newcommand{\norm}[1]{\left|\left|#1\right|\right|}
\newcommand{\A}{\mathcal{A}}
\newcommand{\B}{\mathcal{B}}
\newcommand{\C}{\mathcal{C}}
\newcommand{\D}{\mathcal{D}}
\newcommand{\G}{\mathcal{G}}
\newcommand{\N}{\mathcal{N}}
\newcommand{\R}{\mathcal{R}}
\newcommand{\bdd}{\mathcal{L}}
\newcommand{\CL}{\C_\Lambda}
\newcommand{\diag}{\operatorname{diag}}
\newcommand{\rp}{\operatorname{Re}}

\begin{document}

\title[Robust Regulation of Port-Hamiltonian Systems]{Robust Regulation of Infinite-Dimensional Port-Hamiltonian Systems}

\thispagestyle{plain}

\author{Jukka-Pekka Humaloja} \author{Lassi Paunonen}
\address{Department of Mathematics, Tampere University of Technology, PO.\ Box 553, 33101 Tampere, Finland}
\email{jukka-pekka.humaloja@tut.fi} \email{lassi.paunonen@tut.fi}

\begin{abstract}
We will give general sufficient conditions under which a  controller achieves robust regulation for a boundary control and observation system. Utilizing these conditions we construct a minimal order robust controller for an arbitrary order impedance passive linear port-Hamiltonian system. The theoretical results are illustrated with a numerical example where we implement a controller for a one-dimensional Euler-Bernoulli beam with boundary controls and boundary observations. 
\end{abstract}

\thanks{Lassi Paunonen is funded by the Academy of Finland grant number 298182.}

\maketitle

\section{Introduction} \label{sec:introduction}

The class of infinite-dimensional port-Hamiltonian systems includes models of flexible systems, traveling waves, heat exchangers, bioreactors, and in general, lossless and dissipative hyperbolic systems on one-dimensional spatial domains \cite{ Aug_Phd, LeGZwa05, JacZwa_Book}. In this paper, we consider robust output regulation for port-Hamiltonian systems in general, and as an example we implement a robust controller for Euler-Bernoulli beam which can be formulated as a second-order port-Hamiltonian system. By robust regulation we mean that the controller asymptotically tracks the reference signal, rejects the disturbance signal and allows some perturbations in the plant.

The {\itshape internal model principle} is the key to understanding how control systems can be robust, i.e., tolerate perturbations in the parameters of the system. The principle indicates that the regulation problem can be solved by including in the controller a suitable {\itshape internal model} of the dynamics of the exosystem that generates the reference and disturbance signals. One of the first robust controllers that utilize the internal model principle is the low-gain controller proposed by Davison \cite{Dav_CDC75}. Davison's controller has many practical advantages as it has simple structure and it can be tuned with input-output measurements. The controller was generalized to infinite-dimensional systems and its tuning process was simplified in \cite{HamPoh00, HamPoh02}.

The main contribution of this paper is that we present sufficient criteria for a controller to achieve robust output regulation for boundary control and observation systems. A corresponding result has already been shown for various system classes \cite{HamPoh10, PauPoh10,  PauPoh14a} but not for boundary control systems. As our second main result, we will construct a minimal order robust regulating controller for an arbitrary order impedance passive linear port-Hamiltonian system for which we can show certain assumptions to hold.

Robust output regulation of port-Hamiltonian systems has been considered by the authors in \cite{HumPau_ECC16, HumPau_MTNS16} where first- and even-order port-Hamiltonian systems were considered, respectively. Outside robust regulation, stability, stabilization and dynamic boundary control of port-Hamiltonian systems have been considered, e.g., in \cite{AugJac14, MacLeG17, RamLeG14, VilZwa09}. This paper generalizes the results of \cite{HumPau_ECC16, HumPau_MTNS16} for port-Hamiltonian systems of arbitrary order $N$. Furthermore, as opposed to \cite{HumPau_ECC16, HumPau_MTNS16} considering only impedance energy preserving systems, here we will be considering impedance passive systems as well. Additionally, here the observation operator is allowed to be unbounded, which is essential for true boundary observation. This is also an extension to the results of \cite{HamPoh02} where robust regulation of boundary control systems with bounded observations was considered.

Robust regulation has been considered for boundary control systems in \cite{HamPoh02} and for well-posed systems in general in \cite{RebWei03}. In both references, the robust regulation result is formulated for a single controller structure, whereas our result (Theorem \ref{thm:gcondtororp}) holds for any controller that includes a suitable internal model of the exosystem and stabilizes the closed-loop system. Furthermore, both references assume that the controlled system is initially stable, which is not required here. We also note that in the proof of Theorem \ref{thm:rrc} we could utilize the frequency domain proof of \cite[Thm. 1.1]{RebWei03} to show that the minimal order controller stabilizes the closed-loop system, but we present an alternative time domain proof instead.

The structure of the paper is as follows. In Section \ref{sec:pec} we present the control system consisting of the plant, exosystem and controller. In Section \ref{sec:rorp} we formulate the robust output regulation problem and present the robust regulation result for boundary control and observation systems. In Section \ref{sec:phs} we present the specific structure of linear port-Hamiltonian systems with stability and stabilization results, so that in Section \ref{sec:rrc} we can construct a robust regulating controller - that is also of minimal order - for these systems. The theoretical results are illustrated in Section \ref{sec:example} where we construct a robust regulating controller for Euler-Bernoulli beam.

Here $\bdd(X,Y)$ denotes the set of bounded linear operators from the normed space $X$ to the normed space $Y$. The domain, range, kernel, spectrum and resolvent of a linear operator $A$ are denoted by $\D(A),\R(A),\N(A),\sigma(A)$ and $\rho(A)$, respectively. The resolvent operator is given by $R(\lambda,A) = (\lambda - A)^{-1}$, and it exists for all $\lambda \in \rho(A)$. The growth bound of the $C_0$-semigroup $T_A(t)$ generated by $A$ is denoted by $\omega_0(T_A)$, and $T_A$ is exponentially stable if $\omega_0(T_A) < 0$. In that case we also say that $A$ is exponentially stable.

\section{The plant, exosystem and controller} \label{sec:pec}

The plant is a boundary control system of the form
\begin{subequations}
\begin{align}
 \dot{x}(t) & = \A x(t), \quad x(0) = x_0, \label{eq:Ap} \\
	\B x(t) & = u(t) + w(t), \label{eq:up} \\
  \C x(t) & = y(t) \label{eq:yp}
\end{align}
	\label{eq:plant}%
\end{subequations}
where the disturbance signal $w(t)$ is generated by the exosystem that will be presented shortly.  In Section \ref{sec:phs}, we will make an additional assumption that the plant is an impedance passive port-Hamiltonian system, but for now it is sufficient to consider the plant a {\itshape boundary control and observation system} given by the following definition:
\begin{definition} \textup{\cite[Def. 2.3.13]{Aug_Phd}}
Let $X,U$ and $Y$ be Hilbert spaces. The system $(\A, \B, \C)$ of linear operators $\A: \D(\A) \subset X \to X$, $\B: \D(\B) \subset X \to U$ and $\C: \D(\C) \subset X \to Y$ is called a boundary control and observation system if the following hold:
\begin{enumerate}
	\item $\D(\A) \subset \D(\B)$ and $\D(\A) \subset \D(\C)$.
	\item The restriction $A = \A|_{\N(\B)}$ of $\A$ to the kernel of $\B$ generates a $C_0$-semigroup $(T_A(t))_{t \geq 0}$ on $X$.
	\item There is a right inverse $B \in \bdd(U,X)$ of $\B$ such that $\R(B) \subset \D(\A)$, $\A B \in \bdd(U,X)$ and $\B B = I_U$.
	\item The operator $\C$ is bounded from $\D(A)$ to $Y$, where $\D(A)$ is equipped with the graph norm of $A$.
\end{enumerate}
\label{def:bcos}
\end{definition}

Let $A = \A|_{\N(\B)}$ be the generator of a $C_0$-semigroup $T_A(t)$ on $X$. 
We define the $\Lambda$-extension $\CL$ of $\C$ by
$$
	\CL x = \lim_{\lambda \to \infty} \lambda \C R(\lambda, A)x,
$$
and its domain $\D(\CL)$ consists of those $x \in X$ for which the limit exists. Throughout this paper, we also assume that $\C$ is {\itshape admissible} for $A$ \cite[Def. 4.3.1]{TucWei_Book}, i.e., for some $\tau > 0$ there exists a constant $K_\tau$ such that
$$
\int_0^\tau ||\C T_A(t)x_0||_Y^2dt \leq K_\tau^2||x_0||_X^2 \quad \forall x_0 \in \D(A).
$$
Furthermore, if there exists a constant $K > 0$ such that $K_\tau \leq K$ for all $\tau > 0$, then we say that $\C$ is {\itshape infinite-time admissible} for $A$, for which we will give sufficient conditions in the port-Hamiltonian context later on.

The exosystem that generates the boundary disturbance signal $w(t)$ and the reference signal $y_{ref}(t)$ is a linear system
\begin{subequations}
\begin{align}
 \dot{v}(t) & = Sv(t), \qquad v(0) = v_0, \\
 w(t) & = Ev(t), \\
 y_{ref}(t) & = -Fv(t)
\end{align}
\label{eq:ex}%
\end{subequations}
on a finite-dimensional space $W = \mathbb{C}^q$ with some $q \in \mathbb{N}$. Here $S \in \bdd(W) = \mathbb{C}^{q \times q}$, $E \in \bdd (W, U)$ and $F \in \bdd (W,Y)$. Furthermore, we assume that $S$ has purely imaginary eigenvalues $\sigma(S) = \{i\omega_k\}_{k=1}^q \subset i\mathbb{R}$ with algebraic multiplicity one.

The transfer function of the plant \eqref{eq:plant} is given by
\begin{equation}
	P(\lambda) = \CL R(\lambda,A)(\A B - \lambda B) + \CL B \in \bdd(U,Y),
\label{eq:P}
\end{equation}
and it is defined for every $\lambda \in \rho(A)$ as $\R(B) \subset \D(\A) \subset \D(\C)$. Note that the boundedness of the transfer function implies that $\lambda \hat{u}$ must be bounded for every $\lambda \in \rho(A)$. Hence, by the Plancherel theorem we must have $u \in H^1$, which we will show to hold at the end of this section. Furthermore, we need to assume that $P(i\omega_k)$ is surjective for all $k \in \{1,2,\ldots,q\}$, which is crucial to the solvability of the robust output regulation problem presented in Section \ref{sec:rorp}. Note that the surjectivity assumption implies that we must have $\dim(U) \geq \dim(Y)$.

Since the plant is a boundary control and observation system, it follows from Definition \ref{def:bcos} that we can define an operator $G := BE \in \bdd(W,X)$ satisfying $\A G \in \bdd(W,X)$, $\B G = E$ and $\R(G) \subset \D(\C)$. It is easily seen by following the proof of \cite[Thm. 3.3.3]{CurZwa_Book} that if $u \in C^2(0, \tau; U)$ and $v \in C^2(0, \tau; W)$ for all $\tau > 0$, then the abstract differential equation
\begin{equation}
	\dot{\xi}(t) = A\xi(t) + \A Bu(t) - B\dot{u}(t) + \A Gv(t) - G \dot{v}(t)
\label{eq:abcs}
\end{equation}
with $\xi(0) = \xi_0$ is well-posed. Furthermore, if $\xi_0 = x_0 - Bu_0 - Gv_0 \in \D(A)$, the classical solutions of \eqref{eq:plant} and \eqref{eq:abcs} are related by $\xi(t) = x(t) - Bu(t) - Gv(t)$, and they are unique.

The plant \eqref{eq:plant} - as well as equation \eqref{eq:abcs} - has a well-defined mild solution for $\dot{u} \in L^p(0, \tau; U)$, $\dot{v} \in L^p(0, \tau ; W)$ for some $p \geq 1$ and $x_0 \in X$. In that case, the summary related to \cite[Thm. 3.3.4]{CurZwa_Book} implies that the mild solution of \eqref{eq:plant} is given by
$$
\begin{aligned}
	  x(t) = & \; T_A(t)(x_0 - Bu_0 - Gv_0) + Bu(t) + Gv(t) \, + \\ & \int_0^t T_A(t-s)(\A Bu(s) - B \dot{u}(s) + \A Gv(t) - G \dot{v}(t))ds.
	  \end{aligned}
$$
Similarly for every $\xi_0 = x_0 - Bu_0 - Gv_0 \in X$, one obtains the mild solution of \eqref{eq:abcs} using the above solution and the relation between $x(t)$ and $\xi(t)$. We will show at the end of this section that $\dot{u} \in L^2(0, \tau; U)$, which together with the fact that $v \in C^\infty(0, \tau; W)$ ensures that the mild solutions are well-defined.

The dynamic error feedback controller is of the form
\begin{subequations}
	\begin{align}
	\dot{z}(t) & = \mathcal{G}_1 z(t) + \mathcal{G}_2 (y(t) - y_{ref}(t)), \qquad z(0) = z_0, \label{eq:z} \\
	 u(t) & = Kz(t) \label{eq:controller_u}
	\end{align}
\label{eq:controller}%
\end{subequations}
on a Banach space $Z$. The parameters $\mathcal{G}_1 \in \bdd(Z)$, $\mathcal{G}_2 \in \bdd(Y,Z)$ and $K \in \bdd(Z,U)$ are to be chosen such that robust output regulation is achieved for the plant \eqref{eq:plant}.

We are finally in the position to give the formulation of the closed-loop system consisting of the plant \eqref{eq:plant} written as the abstract differential equation \eqref{eq:abcs} and the controller \eqref{eq:controller}. Furthermore, we will show that $\dot{u} \in L^2(0,\tau;U)$ for every $\tau > 0$. Using the above notation and definitions, the closed-loop system can be written on the extended state space $X_e = X \times Z$ with the extended state $\xi_e(t) = (\xi(t), z(t))^T$ as
\begin{subequations}
\begin{align}
 & \dot{\xi}_e(t) = A_e\xi_e(t) + B_e v(t), \quad \xi_e(0) = \xi_{e0},
  \label{eq:extended_state} \\
 & e(t) = C_e \xi_e(t) + D_e v(t), \label{eq:extended_output}
\end{align}
\label{eq:extended}%
\end{subequations}
where $e(t) := y(t) -y_{ref}(t)$ is the regulation error, $\xi_{e0} = (\xi_0, z_0)^T$, $C_e = [\CL \quad \CL BK]$, $D_e = \CL G + F$ and
$$
\begin{aligned}
& A_e = \left[ \begin{array}{cc} A - BK\G_2 \CL & \A BK - BK(\G_1 + \G_2 \CL BK) \\
	\G_2 \CL & \G_1 + \G_2 \CL BK
\end{array} \right], \\
& B_e = \left[ \begin{array}{c} \A G - GS - BK\G_2 (\CL G + F) \\ \G_2(\CL G + F)
 \end{array} \right].
\end{aligned}
$$

The operator $A_e$ has domain $\D(A_e) = \D(A) \times Z$, and it can be written in the form
$$
\begin{aligned}
	A_e = & \left[ \begin{array}{cc} A & 0 \\ 0 & \G_1 \end{array} \right] + \left[ \begin{array}{c} -BK\G_2 \\
	\G_2 \end{array} \right]\left[\CL \quad \CL BK\right] + \left[ \begin{array}{cc} 0 & \A BK - BK\G_1 \\ 0 & 0 \end{array} \right] \\
	:= & A_1 + A_2C_e + A_3.
\end{aligned}
$$
Since all the operators associated with the controller \eqref{eq:controller} are bounded and since $\A B,B \in \bdd(U,X)$ due to the plant \eqref{eq:plant} being a boundary control and observation system, it follows that the operators $A_2$ and $A_3$ are bounded. Furthermore, since $\C$ is admissible for $A$ and $\CL B \in \bdd(U,Y)$, it follows that $C_e$ is admissible for $A_1$. Thus, since $A_1$ is clearly the generator of a $C_0$-semigroup, $A_2$ and $A_3$ are bounded, and $C_e$ is admissible for $A_1$, it follows from \cite[Thm. 5.4.2]{TucWei_Book} and standard perturbation theory that the operator $A_e$ is the generator of a $C_0$-semigroup, and that $C_e$ is admissible for $A_e$ as well. Finally, combining \eqref{eq:controller} and \eqref{eq:extended_output} we obtain that $\dot{u} = K\G_1 z + K\G_2 (C_e \xi_e + D_e v)$, which by the above reasoning shows that $\dot{u} \in L^2(0,\tau; U)$ for all $\tau > 0$, and thus, the mild solutions of \eqref{eq:plant} and \eqref{eq:abcs} are well-defined.

\section{The robust output regulation problem and the internal model principle} \label{sec:rorp}

In this section, we formulate the {\itshape robust output regulation problem} and present the concept of the internal model via the {\itshape $\G$-conditions}. After that, we are in the position to present and prove the first main result of this paper.

In order to discuss robustness, we consider perturbations $(\tilde{\A}, \tilde{\B}, \tilde{\C}, \tilde{E}, \tilde{F}) \in \mathcal{O}$ of the operators $(\A, \B, \C, E, F)$. The class $\mathcal{O}$ of perturbations is defined such that the perturbed operators $(\tilde{\A}, \tilde{\B}, \tilde{\C}, \tilde{E}, \tilde{F})$ satisfy the following assumptions which the operators $(\A, \B, \C, E, F)$ are assumed to satisfy as well.
\begin{assumption}
The operators $(\tilde{\A},\tilde{\B},\tilde{\C},\tilde{E},\tilde{F})$ satisfy the following:
\begin{enumerate}
	\item The plant $(\tilde{\A},\tilde{\B},\tilde{\C})$ is a boundary control and observation system.
	\item The operator $\tilde{\C}$ is admissible for $\tilde{A} = \tilde{\A}|_{\N(\tilde{\B})}$.
	\item The transfer function of the plant $(\tilde{\A},\tilde{\B},\tilde{\C})$ is surjective and bounded for every eigenvalue of $S$.
	\item $\tilde{E} \in \bdd(W,U)$ and $\tilde{F} \in \bdd(W,Y)$.
\end{enumerate}
\label{ass:standing}
\end{assumption}
\noindent It is easy to see that these conditions are satisfied for arbitrary bounded perturbations to $E$ and $F$, whereas the boundary control and observation system requirement imposes stricter conditions on the perturbations on $\A,\B$ and $\C$. However, at least sufficiently small bounded perturbations are acceptable. Note that the operators $B$ and $G$ associated with the boundary control and observation system will also change when the system is perturbed. We denote these operators by $\tilde{B}$ and $\tilde{G}$.

\noindent {\bf The Robust Output Regulation Problem.} Choose a controller $(\mathcal{G}_1, \mathcal{G}_2, K)$ in such a way that the following are satisfied:
\begin{enumerate}
	\item The closed-loop system generated by $A_e$ is exponentially stable.
	\item For all initial states $\xi_{e0} \in X_e$ and $v_0 \in W$, the regulation error satisfies $e^{\alpha \cdot} e(\cdot) \in L^2(0,\infty; Y)$ for some $\alpha > 0$.
	\item If $(\A, \B, \C, E, F)$ are perturbed to $(\tilde{\A}, \tilde{\B}, \tilde{\C}, \tilde{E}, \tilde{F}) \in \mathcal{O}$ in such a way that the closed-loop system remains exponentially stable, then for all initial states $\xi_{e0} \in X_e$ and $v_0 \in W$, the regulation error satisfies $e^{\tilde{\alpha}\cdot} e(\cdot) \in L^2(0,\infty; Y)$ for some $\tilde{\alpha} > 0$.
\end{enumerate}
\noindent We note that without the last item in the above list the problem is called {\itshape output regulation problem} which will be considered in the proof of our first main result in the next subsection.

The internal model principle states that the robust output regulation problem can be solved by including a suitable internal model of the dynamics of the exosystem in the controller. The internal model can be characterized using the definition of $\G$-conditions below. What follows is our first main result where we show that a controller satisfying the $\G$-conditions is robust.
\begin{definition} \textup{\cite[Def. 10]{HamPoh10}}
	A controller $(\G_1, \G_2, K)$ is said to satisfy the $\G$-conditions if
	\begin{subequations}
		\begin{align}
		\R(i\omega_k - \G_1) \cap \R(\G_2) & = \{0\}, \\
		\N(\G_2) & = \{0\}
		\end{align}
\label{eq:gconditions}%
\end{subequations}
for all $k \in \{1,2,\ldots,q\}$, where $\sigma(S) = \{i\omega_k\}_{k=1}^q$.
	\label{def:gcond}
\end{definition}

\subsection{Sufficient Robustness Criterion for a Controller} \label{rorp:thm}

We will now show that a controller $(\G_1, \G_2, K)$ that satisfies the $\G$-conditions solves the robust output regulation problem for a boundary control and observation system, provided that the controller exponentially stabilizes the closed-loop system.

\begin{theorem}
Assume that a controller $(\G_1, \G_2, K)$ exponentially stabilizes the closed-loop system. If the controller satisfies the $\G$-conditions, then it solves the robust output regulation problem. The controller is guaranteed to be robust with respect to all perturbations under which the closed-loop system remains exponentially stable and Assumption \ref{ass:standing} is satisfied.

\begin{proof}
Let $(\tilde{\A}, \tilde{\B}, \tilde{\C}, \tilde{E}, \tilde{F})$ be arbitrary perturbations of class $\mathcal{O}$ such that the perturbed closed-loop system generated by $\tilde{A}_e$ is exponentially stable. As the perturbations of the class $\mathcal{O}$ satisfy Assumption \ref{ass:standing}, it follows that $\tilde{B}_e$ and $\tilde{D}_e$ are bounded and $\tilde{C}_e$ is admissible for $\tilde{A}_e$. Thus, the closed-loop system is a regular linear system, and by \cite[Thm. 4.1]{PauPoh14a} we have that the controller $(\G_1,\G_2,K)$ solves the output regulation problem if and only if the regulator equations $\Sigma S = \tilde{A}_e\Sigma + \tilde{B}_e$ and $0 = \tilde{C}_e\Sigma + \tilde{D}_e$ have a solution $\Sigma := (\Pi, \Gamma)^T \in \bdd(W,X_e)$. Note that the result of \cite[Thm. 4.1]{PauPoh14a} only requires that the closed-loop system is regular, and therefore it can be used here. Further note that as $\tilde{A}_e$ is assumed to be exponentially stable and $\sigma(S) \subset i\mathbb{R}$, by \cite{Pho91} the Sylvester equation $\Sigma S = \tilde{A}_e \Sigma + \tilde{B}_e$ has a unique solution $\Sigma \in \bdd(W, X_e)$ satisfying $\R(\Sigma) \subset \D(\tilde{A}_e)$. Thus, in order to show that the controller solves the output regulation problem, it remains to show that the bounded solution $\Sigma$ of the Sylvester equation satisfies the second regulator equation as well. We will do this for the arbitrary perturbations $(\tilde{\A}, \tilde{\B}, \tilde{\C}, \tilde{E}, \tilde{F}) \in \mathcal{O}$, which implies that the controller is robust under these perturbations, i.e., it solves the robust output regulation problem.

Let $k \in \{1,2,\ldots,q\}$ be arbitrary and consider the eigenvector $\phi_k$ of $S$ associated with the corresponding eigenvalue $i\omega_k$  satisfying $S\phi_k = i\omega_k\phi_k$. Then $\Sigma S \phi_k = \tilde{A}_e\Sigma\phi_k + \tilde{B}_e \phi_k$ implies $(i\omega_k - \tilde{A}_e)\Sigma \phi_k = \tilde{B}_e \phi_k$, which yields
$$
\begin{aligned}
& \left[ \begin{array}{c}
	(i\omega_k - \tilde{A} + \tilde{B}K\G_2\tilde{\C}_\Lambda)\Pi\phi_k - (\tilde{\A} \tilde{B}K - \tilde{B}K(\G_1 + \G_2\tilde{\C}_\Lambda\tilde{B}K))\Gamma\phi_k \\
	-\G_2 \tilde{\C}_\Lambda\Pi\phi_k + (i\omega_k - \G_1)\Gamma\phi_k - \G_2 \tilde{\C}_\Lambda\tilde{B}K\Gamma\phi_k
\end{array} \right]\\ 
& = \left[ \begin{array}{c}
	(\tilde{A}\tilde{G} - \tilde{G}S - \tilde{B}K\G_2(\tilde{\C}_\Lambda\tilde{G}+ \tilde{F}))\phi_k \\
	\G_2(\tilde{\C}_\Lambda\tilde{G}+ \tilde{F})\phi_k
\end{array} \right].
\end{aligned}
$$
The second line implies $(i\omega_k - \G_1)\Gamma\phi_k = \G_2(\tilde{\C}_\Lambda\Pi\phi_k + \tilde{\C}_\Lambda\tilde{B}K\Gamma\phi_k + (\tilde{\C}_\Lambda\tilde{G} + \tilde{F})\phi_k)$, and now by the $\G$-conditions we have that $0 = \tilde{\C}_\Lambda\Pi\phi_k + \tilde{\C}_\Lambda\tilde{B}K\Gamma\phi_k + (\tilde{\C}_\Lambda\tilde{G}+ \tilde{F})\phi_k = \tilde{C}_e\Sigma\phi_k + \tilde{D}_e\phi_k$. As the eigenvectors $\phi_k$ form an orthogonal basis on $W$ and the choice of $k$ was arbitrary, it follows that $\Sigma$ satisfies the second regulator equation $\tilde{C}_e \Sigma + \tilde{D}_e = 0$ as well. Thus, \cite[Thm. 4.1]{PauPoh14a} implies that the controller solves the robust output regulation problem.
\end{proof}
\label{thm:gcondtororp}
\end{theorem}

\section{Background on port-Hamiltonian systems} \label{sec:phs}

In this section, we give some background to port-Hamiltonian systems. We note that while \cite{LeGZwa05} is a classical reference paper regarding these systems, we use \cite{Aug_Phd} as our main reference as it gives a slightly more general formulation for port-Hamiltonian systems than \cite{LeGZwa05}. Therefore we will cite \cite{Aug_Phd} for the base results as well, even though essentially the same results can be found in \cite{LeGZwa05}.

Define a linear port-Hamiltonian operator $\A$ of order $N$ on the spatial interval $\zeta \in [a,b]$ as follows:
\begin{definition} \textup{\cite[Def. 3.2.1]{Aug_Phd}}
Let $N \in \mathbb{N}$ and $P_k \in \mathbb{C}^{n \times n}$ satisfying $P^*_k = (-1)^{k+1}P_k$ for $k \in \{1,2\ldots,N\}$ with $P_N$ invertible. Furthermore, let $P_0 \in L^\infty(a,b;\mathbb{C}^{n \times n})$ satisfying $\rp(P_0(\zeta)) := \frac{1}{2}(P_0(\zeta) + P_0^*(\zeta)) \leq 0$ for a.e. $\zeta \in [a, b]$. Let the state space $X = L^2(a,b;\mathbb{C}^n)$ be equipped with the inner product $\langle \cdot, \cdot \rangle_X = \langle \cdot, \mathcal{H}\cdot \rangle_{L^2}$ where $\mathcal{H}: [a,b] \to \mathbb{C}^{n \times n}$ satisfies $m|\xi|^2 \leq \langle \xi, \mathcal{H}(\zeta)\xi\rangle_{\mathbb{C}^n} \leq M|\xi|^2, \quad \xi \in \mathbb{C}^n$  a.e. $\zeta \in [a,b]$ for some constants $0 < m \leq M < \infty$. Then the operator $\A: \D(\A) \subset X \to X$ defined as
$$
	\A x(\zeta, t) :=  \sum_{k=1}^N P_k \frac{\partial^k}{\partial \zeta^k} (\mathcal{H}(\zeta)x(\zeta, t)) + P_0(\zeta) \mathcal{H}(\zeta)x(\zeta,t),
$$
with domain $\D(\A) = \{x \in X: \mathcal{H}x \in H^N(a,b; \mathbb{C}^n)\}$ is called a linear port-Hamiltonian operator of order $N$.
\label{def:pHo}
\end{definition}
\noindent Let $\Phi : H^N(a,b; \mathbb{C}^n) \rightarrow \mathbb{C}^{2nN}$ defined by
$$
	\Phi(x) := (x(b), \ldots, x^{(N-1)}(b),x(a),\ldots, x^{(N-1)}(a))^T
$$
be the boundary trace operator and define the boundary port variables $f_\partial$, $e_\partial$ by
\begin{equation}
	\left[ \begin{array}{c}
	f_\partial \\ e_\partial 
\end{array} \right] := \frac{1}{\sqrt{2}} \left[ \begin{array}{rr}
	Q & -Q \\ I & I 
\end{array} \right] \Phi(\mathcal{H}x) := R_{ext} \Phi(\mathcal{H}x)
\label{eq:bef}
\end{equation}
where $Q \in \mathbb{C}^{nN \times nN}$ is a block matrix given by
$$
	Q_{ij} := \begin{cases} 
		(-1)^{j-1} P_{i+j-1},	& i + j \leq N + 1 \\ 0, & \mbox{else}
	\end{cases}.
$$
Note that since $P_N$ is assumed to be invertible, it follows that $Q$ is invertible, and hence, $R_{ext}$ is invertible as well. 

Using the boundary port variables we can now define the boundary control and boundary observation operators  $\mathcal{B}$ and $\mathcal{C}$, respectively. Their definitions are included in the following definition of port-Hamiltonian systems.

\begin{definition} \textup{\cite[Def. 3.2.10]{Aug_Phd}}
Let $\A$ be a port-Hamiltonian operator of order $N$ with associated boundary port variables $f_\partial$ and $e_\partial$. Further let $W_B,W_C \in \mathbb{C}^{nN \times 2nN}$ be full rank matrices such that $\N(W_B) \cap \N(W_C) = \{0\}$. Then the input map $\B:\D(\B) = \D(A) \subset X \to U := \mathbb{C}^{nN}$ and the output map $\C: \D(\C) = \D(A) \subset X \to Y := \mathbb{C}^{nN}$ are defined as
\begin{equation}
	\mathcal{B}x(t) := W_B \left[ \begin{array}{c} f_\partial(t) \\ e_\partial(t) \end{array} \right], \quad
	\mathcal{C}x(t) := W_C \left[ \begin{array}{c} f_\partial(t) \\ e_\partial(t) \end{array} \right]
	\label{eq:bcdef}
\end{equation} 
and the system $(\A,\B,\C)$ is called a port-Hamiltonian system.
\label{def:pHs}
\end{definition}

\noindent We note that the above definition implies that we have full control and measurements, which is not very common in practice. However, the exponential stability criterion given in part \ref{wbkwc:1}) of Lemma \ref{lem:wbkwc} essentially requires that we have full control. If we were considering a less general class of port-Hamiltonian systems, e.g., first- or even order systems, we could utilize \cite[Thm. III.2]{VilZwa09} or \cite[Prop. 2.16]{AugJac14}, respectively, to obtain exponential stability with fewer controls. However, to our knowledge there are no weaker exponential stability criteria than the one given in part \ref{wbkwc:1}) of Lemma \ref{lem:wbkwc} for arbitrary order port-Hamiltonian systems, and thus, we assume having full control and measurements.

We have by \cite[Thm. 3.2.21]{Aug_Phd} that a port-Hamiltonian system $(\A,\B,\C)$ is a boundary control and observation system if and only if the operator $A = \A|_{\N(\B)}$ generates a $C_0$-semigroup on $X$. Furthermore, by \cite[Thm. 3.3.6]{Aug_Phd} the operator $A$ generates a contractive $C_0$-semigroup if and only if $W_B \Sigma W_B^* \geq 0$ where
\begin{equation}
	\Sigma := \left[ \begin{array}{cc} 0 & I \\ I & 0 \end{array} \right].
\label{eq:Sigma}
\end{equation}
Following \cite[Def. 3.2.12]{Aug_Phd}, we define a system $(\A,\B,\C)$ {\itshape impedance passive} if it satisfies
$$
\rp \langle \A x(t), x(t) \rangle_X \leq  \rp \langle \B x(t) , \C x(t) \rangle_{\mathbb{C}^{nN}}, \quad x \in \D (\A)
$$
and {\itshape impedance energy preserving} if the above holds as an equality.  These systems can be easily identified based on $W_B, W_C$ and $P_0$. Define a matrix $P_{W_B,W_C}$ such that 
$$
	P_{W_B,W_C}^{-1} = \left[ \begin{array}{cc} W_B\Sigma W_B^* & W_B\Sigma W_C^* \\
	W_C\Sigma W_B^* & W_C\Sigma W_C^*
	\end{array} \right].
$$
By \cite[Prop. 3.2.16]{Aug_Phd}, a port-Hamiltonian system is impedance energy preserving if and only if $P_0(\zeta) = -P_0(\zeta)^*$ for a.e. $\zeta \in [a, b]$ and $P_{W_B,W_C} = \Sigma$, and it is impedance passive if and only if $\rp P_0(\zeta) \leq 0$ for a.e. $\zeta \in [a, b]$ and $P_{W_B,W_C} \leq \Sigma$.

We consider impedance energy preserving and impedance passive port-Hamiltonian systems as they can be exponentially stabilized using output feedback. Stabilization of port-Hamiltonian systems with negative output feedback was first presented for first-order impedance energy preserving port-Hamiltonian systems in \cite[Sec. IV]{VilZwa_CDC05}, and we will next generalize the result for systems of arbitrary order $N$.
\begin{lemma} \hfill
\begin{enumerate}[a)]
\item A port-Hamiltonian system that satisfies $W_B\Sigma W_B^* > 0$ is exponentially stable. \label{wbkwc:1}
\item An impedance passive port-Hamiltonian system can be exponentially stabilized using negative output feedback $u(t) = -\kappa y(t)$ for any $\kappa > 0$. \label{wbkwc:2}
\end{enumerate}
\begin{proof} a) The claim can be proved similarly to \cite[Lem. 2]{HumPau_MTNS16} by using the techniques utilized in the proof of \cite[Lem. 9.1.4]{JacZwa_Book} and the estimate $\rp\langle Ax,x\rangle_X \leq \operatorname{Re}\langle f_\partial, e_\partial \rangle_{\mathbb{C}^{nN}}$  which holds as $\operatorname{Re} P_0(\zeta) \leq 0$ a.e. $\zeta \in [a, b]$. Eventually, we obtain 
$$
\operatorname{Re}\langle Ax,x\rangle_X \leq -\gamma \sum_{k=0}^{N-1}\sum_{\zeta = a,b}|(\mathcal{H} x)^{(k)}(\zeta)|^2
$$
for some $\gamma > 0$, which by \cite[Thm. 4.3.24]{Aug_Phd} is sufficient for the port-Hamiltonian system being exponentially stable. 

b) Let $W_B$ and $W_C$ be such that the port-Hamiltonian system is impedance passive. It has been shown in \cite[Sec. IV]{VilZwa_CDC05} that the closed-loop system with negative output feedback $u(t) = -\kappa y(t)$ can be written as
$$
\begin{aligned}
\dot{x}(t) & = \A x(t), \\
(W_B + \kappa W_C) \left[ \begin{array}{c} f_\partial (t) \\ e_\partial (t) 
\end{array} \right] & = (\B + \kappa \C)x (t) \equiv 0, \\
\C x(t) & = y(t).
\end{aligned} 
$$
By \cite[Prop. 3.2.16, Lem. 3.2.18]{Aug_Phd}, it holds for impedance passive port-Hamiltonian systems that $W_B\Sigma W_B^* \geq 0$, $W_C\Sigma W_C^* \geq 0$ and $W_B\Sigma W_C^* = I = W_C\Sigma W_B^*$. Denote $W_\kappa := W_B + \kappa W_C$ which satisfies
$$
	W_\kappa\Sigma W_\kappa^* = W_B\Sigma W_B^* + 2\kappa I + \kappa^2W_C\Sigma W_C^* \geq 2\kappa I > 0,
$$
and now part \ref{wbkwc:1}) completes the proof.
\end{proof}
\label{lem:wbkwc}
\end{lemma}

\section{Robust regulating controller for impedance passive port-Hamiltonian systems} \label{sec:rrc}

In this section, we will construct a finite dimensional, minimal order controller for an impedance passive port-Hamiltonian system and a finite dimensional exosystem as given in \eqref{eq:ex}. The choices of the controller parameters $(\G_1, \G_2, K)$ are adopted from \cite[Sec. 4]{Pau16}. However, as an impedance passive port-Hamiltonian system is not necessarily exponentially stable to begin with, we will need to add an extra term to the controller in order to ensure the exponential stability of the closed-loop system. The controller that we will construct is of the form
$$
\begin{aligned}
	\dot{z}(t) & = \G_1z(t) + \G_2e(t), \quad z(0) = z_0, \\
	u(t) & = Kz(t) - \kappa e(t),
\end{aligned}
$$
where as opposed to the controller given in \eqref{eq:controller} we have the extra feedthrough term $-\kappa e(t)$. Here the control signal consists of two parts $u(t) = u_1(t) + u_2(t)$ where the second term contributes to exponentially stabilizing the plant and the first one provides the robust regulating control. Note that instead of $-\kappa y(t)$ we use $-\kappa e(t)$ which we will show to stabilize the plant as well. Furthermore, using $-\kappa e(t)$ simplifies the controller as $y(t)$ and $y_{ref}(t)$ are not needed separately.

We define $Z = Y^q$. The controller parameters are chosen as $\kappa > 0$ and
$$
\begin{aligned}
	\G_1 & = \diag\left(i\omega_1I_Y,i\omega_2I_Y, \ldots, i\omega_qI_Y\right) \in \bdd(Z), \\
	K & = \epsilon K_0 = \epsilon\left[K_0^1,K_0^2, \ldots, K_0^q\right] \in \bdd(Z,U), \\
	\G_2 & = (\G_2^k)_{k=1}^q = (-(P_\kappa(i\omega_k)K_0^k)^*)_{k=1}^q 
	 \in \bdd(Y,Z)
\end{aligned}
$$
where $\epsilon > 0$ is the tuning parameter and $P_\kappa(i\omega_k) = P(i\omega_k)(I + \kappa P(i\omega_k))^{-1}$ is the transfer function of the triplet $(\A,\B + \kappa\C, \C)$. Note that since $P(i\omega_k)$ is assumed to be surjective for every $k \in \{1,2,\ldots,q\}$, $P_\kappa(i\omega_k)$ is surjective as well. Further note that if we choose $K_0^k = P_\kappa(i\omega_k)^{\dagger}$ (the Moore-Penrose pseudoinverse of $P_\kappa(i\omega_k)$), then $\G_2^k = -I_Y$ for all $k \in \{1,2,\ldots q\}$.

\begin{theorem}
	Assume that $(\A,\B,\C)$ is an impedance passive port-Hamiltonian system of an arbitrary order $N$ and $(S,E,F)$ is a finite-dimensional exosystem such that Assumption \ref{ass:standing} is satisfied. Then there exists an $\epsilon_\kappa^* > 0$ such that for any $0 < \epsilon < \epsilon_\kappa^*$ the controller with the above parameter choices solves the robust output regulation problem.
\begin{proof}
Consider an input of the form $u(t) = Kz(t) - \kappa e(t) = u_1(t) - \kappa y(t) + \kappa y_{ref}(t)$. The plant with such an input can be written as
	$$
	\begin{aligned}
	\dot{x}(t) & = \A x(t), \\
	 (\B + \kappa \C)x(t) & = u_1(t) + \kappa y_{ref}(t) + w(t), \\
	\C x(t) & = y(t),
	\end{aligned}
	$$
where we also included the boundary disturbance signal $w(t)$. Note that as $w(t) = Ev(t)$ and $y_{ref}(t) = -Fv(t)$, the term $\kappa y_{ref}(t)$ can be considered an additional disturbance to the original system. 

We know by Lemma \ref{lem:wbkwc} that the negative output feedback exponentially stabilizes the impedance passive port-Hamiltonian system, and thus, the operator $A_\kappa := \A|_{\N(\B + \kappa \C)}$ generates an exponentially stable $C_0$-semigroup on $X$. Furthermore, as the stabilized plant is a boundary control and observation system, there exists an operator $B_\kappa$ satisfying $(\B + \kappa \C)B_\kappa = I_U$, and we can define an operator $\G_\kappa := B_\kappa(E - \kappa F)$ that satisfies $(\B + \kappa \C)G_\kappa = E - \kappa F$ and takes the reference signal $\kappa y_{ref}(t)$ into account.

The closed-loop system consisting of the plant and the controller is still given as in \eqref{eq:extended} with $A,B$ and $G$ replaced by $A_\kappa, B_\kappa$ and $G_\kappa$, respectively, and the $\Lambda$-extension of $\C$ is given by $\CL x= \lim_{\lambda \to \infty} \lambda \C R(\lambda, A_\kappa)x$. Note that since the plant is an impedance passive port-Hamiltonian system, we have by Lemma \ref{lem:wbkwc} that $(W_B + \kappa W_C)\Sigma (W_B + \kappa W_C)^* > 0$, and thus, by Lemma \ref{lem:admis} presented in the Appendix the operator $\C$ is admissible for $A_\kappa$.

Now that the feedthrough term of the controller is associated with the plant, the remaining controller is of the standard form given in \eqref{eq:controller}. Thus, we have by the proof of \cite[Thm 4.1]{Pau16} that the controller satisfies the $\G$-conditions, and hence, by Theorem \ref{thm:gcondtororp} the controller solves the robust output regulation problem, provided that the closed-loop system is exponentially stable.

To conclude the proof, we will show that the closed-loop operator $A_e$ is similar to an exponentially stable operator and hence, exponentially stable. Choose a similarity transformation
$$
	Q = \left[ \begin{array}{rr} -I & \epsilon H \\ 0 & I 
\end{array} \right] = Q^{-1} \in \bdd(X_e)
$$
where the operator $H := (H_1,H_2,\ldots,H_q) \in \bdd(Z,\D(A_\kappa))$ is chosen as
$$
	H_k := R(i\omega_k,A_\kappa)(\A B_\kappa - i\omega_k B_\kappa)K_0^k 
$$
for all $k \in \{1,2,\ldots,q\}$. Let us define $\hat{A}_e := QA_e Q^{-1}$. We will next show that $\hat{A}_e$ is exponentially stable, which implies that $A_e$ is exponentially stable as well.

By the choices of $H_k$, we have $(i\omega_k - A_\kappa)H_k = \A B_\kappa K_0^k - i\omega_k B_\kappa K_0^k$, i.e., $H_k i\omega_k = A_\kappa H_k + \A B_\kappa K_0^k - B_\kappa K_0^k i\omega_k$, and thus, $H\G_1 = A_\kappa H + \A B_\kappa K_0 - B_\kappa K_0\G_1$ due to the diagonal structure of $\G_1$. Furthermore,
$$
\CL(H_k + B_\kappa K_0^k) = \CL R(i\omega_k,A_\kappa)(\A B_\kappa - i\omega_k B_\kappa)K_0^k + \CL B_\kappa K_0^k = P_\kappa(i\omega_k)K_0^k,
$$ 
and thus, $\CL(H + B_\kappa K_0) = -\G_2^*$. Using the above identities $\hat{A}_e$ can be written as
$$
\resizebox{.99\hsize}{!}{$
\hat{A}_e = \left[ \begin{array}{cc}
	A_\kappa  - \epsilon (H + B_\kappa K_0)\G_2 \CL & 0 \\
	-\G_2 \CL & \G_1 - \epsilon \G_2 \G_2^*
\end{array} \right] + \epsilon^2 \left[ \begin{array}{cc}
	0 & -(H+B_\kappa K_0)\G_2 \G_2^* \\ 0 & 0 
\end{array} \right]. $}
$$

Since $\C$ is admissible for $A_\kappa$ and $(H + B_\kappa K_0)\G_2$ is bounded, there exists an $\epsilon_\kappa > 0$ such that for all $0 < \epsilon < \epsilon_\kappa$ the operator $A_\kappa - \epsilon (H + B_\kappa K_0)\G_2\CL$ generates an exponentially stable semigroup. Furthermore, we have by \cite[App. B]{HamPoh11} that the semigroup generated by $\G_1 - \epsilon \G_2\G_2^*$ is exponentially stable for every $\epsilon > 0$ and that there exists a constant $M > 0$ such that $\norm{R(\lambda, \G_1 - \epsilon \G_2\G_2^*)} \leq M/\epsilon$ for $\lambda \in \mathbb{C}_+$. Consider the operator $\hat{A}_e$ in the form $A_1 + \epsilon^2 A_2$. Using the above upper bound for $\norm{R(\lambda, \G_1 - \epsilon\G_2\G_2^*)}$ it can be shown that there exists an $\epsilon^*$ such that for all $0 < \epsilon < \epsilon^*$ and $\lambda \in \mathbb{C}_+$ we have $\norm{\epsilon^2A_2R(\lambda, A_1)} < 1$. Thus, it follows that there exists an $\epsilon_\kappa^*$ such that for all $0 < \epsilon < \epsilon_\kappa^*$ the resolvent of $\hat{A}_e$ is bounded in the right half plane, i.e., $\hat{A}_e$ is exponentially stable.

Since the controller satisfies the $\G$-conditions and the closed-loop system is exponentially stable for every $0 < \epsilon < \epsilon_\kappa^*$, we have by Theorem \ref{thm:gcondtororp} that the controller solves the robust output regulation problem for any $0 < \epsilon < \epsilon_\kappa^*$.
\end{proof}
\label{thm:rrc}
\end{theorem}

\section{Robust control of a 1D Euler-Bernoulli beam} \label{sec:example}

In this section, we construct a robust controller for Euler-Bernoulli beam which is an example of a port-Hamiltonian system of order two. The formulation of Euler-Bernoulli beam as a port-Hamiltonian system is adopted from \cite[Ex. 3.1.6]{Aug_Phd}.

The Euler-Bernoulli beam equation is given on the spatial interval $\zeta \in [0, 1]$ by
$$
	\rho(\zeta) \frac{\partial^2}{\partial t^2}\nu(\zeta, t) = - \frac{\partial^2}{\partial \zeta^2} \big( EI(\zeta) \frac{\partial^2}{\partial \zeta^2} \nu(\zeta, t) \big)
$$
where $\nu(\zeta, t)$ denotes the displacement at position $\zeta$ at time $t$, $\rho(\zeta)$ is the mass density times the cross sectional area, $E(\zeta)$ is the modulus of elasticity and $I(\zeta)$ is the area moment of the cross section. Due to their physical interpretations, the functions $\rho, E$ and $I$ are uniformly bounded and strictly positive for all $\zeta \in [0, 1]$.

In order to write the Euler-Bernoulli beam equation as a port-Hamiltonian system, let us define the state $x(\zeta, t)$ by
$$
	x(\zeta, t) = \left[ \begin{array}{c}
		x_1(\zeta, t) \\ x_2(\zeta, t) 
	\end{array} \right] := \left[ \begin{array}{r}
			\rho \nu_t(\zeta, t) \\ \nu_{\zeta\zeta}(\zeta, t)
	\end{array} \right].
$$
Now we can write the equation as $\partial_t x(\zeta, t) = \A x(\zeta, t)$ where
$$
	\displaystyle \A x(\zeta, t) := \left[ \begin{array}{rr}
	 0 & -1 \\ 1 & 0
	\end{array} \right]  \frac{\partial^2}{\partial \zeta^2} \left( \left[ \begin{array}{cc}
		\rho(\zeta)^{-1 } & 0 \\ 0 &  EI(\zeta)
	\end{array} \right] x(\zeta, t) \right),
$$
which is a second-order port-Hamiltonian operator with $P_0 = P_1 = 0$,
$$
	P_2 = \left[ \begin{array}{rr}
	0 & -1 \\ 1 & 0
	\end{array} \right] \quad \mbox{and} \quad \mathcal{H}(\zeta) = \left[ \begin{array}{cc}
	\rho(\zeta)^{-1} & 0 \\ 0 & EI(\zeta) 
	\end{array} \right].
$$
Using the new state variables, define the control and observation operators by 
$$
\B x(\cdot, t) := \left[ \begin{array}{c} x_1'(0, t) \\ x_1(0,t) \\ x_2'(1,t) \\ x_2(1,t) \end{array}\right] \quad \mbox{and} \quad \C x(\cdot, t) := \left[ \begin{array}{r}
-x_2(0,t) \\ x_2'(0,t) \\ -x_1(1,t) \\ x_1'(1,t) \end{array}\right], 
$$
from which it can be seen that the triple $(\A, \B, \C)$ is an impedance energy preserving port-Hamiltonian system.

Let the reference signal $y_{ref}$ and the disturbance signal $d$ be given by 
$$
y_{ref}(t) := \left[ \begin{array}{c} - \sin(\pi t) \\ - \cos(2\pi t) \\ \cos(\pi t) \\ \sin(2\pi t) \end{array} \right] \quad \mbox{and} \quad d(t) := \left[ \begin{array}{c}
	\sin(2\pi t) \\ \cos(\pi t) \\ \cos( 2\pi t) \\ \sin(\pi t)
\end{array} \right],
$$
so that we have $S := \operatorname{diag}(-2i\pi, -i\pi, i\pi, 2i\pi)$, and $E$ and $F$ are suitably chosen matrices.

The controller parameters $(\G_1, \G_2, K, \kappa)$ are chosen according to the previous section, i.e., we choose
$$
\begin{aligned}
	\kappa & = 1,  \quad \epsilon = 0.17, \\
	\G_1 & = \diag\left(-2i\pi I_Y,-i\pi I_Y, i\pi I_Y, 2i\pi I_Y\right), \\
	\G_2 & = (\G_2^k)_{k=1}^4, \quad \G_2^k = -I_Y \quad \forall k \in \{1,2,3,4\}, \\
	 K & = \epsilon\left[P_\kappa(-2i\pi)^{-1}, P_\kappa(-i\pi)^{-1}, P_\kappa(i\pi)^{-1}, P_\kappa(2i\pi)^{-1} \right],
\end{aligned}
$$
where $P_\kappa$ is the transfer function of the triplet $(\A, \B + \kappa\C, \C)$ and $\epsilon$ is chosen such that the growth bound of the closed-loop system is close to its minimum. Note that as we chose $K_0^k = P_\kappa(i\omega_k)^{-1}$, each block of $\G_2$ is equal to $-I_Y$.

\begin{figure}[!ht]
\includegraphics[width=\columnwidth]{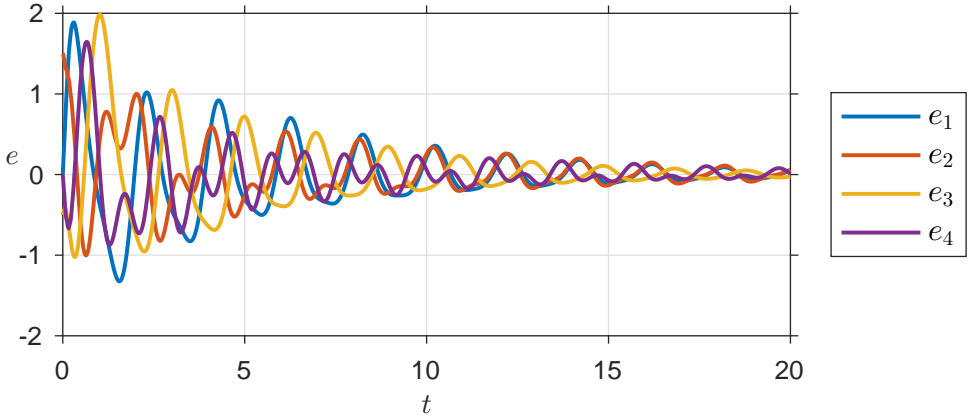}
\caption{Regulation error on $t \in [0, 20]$.}
\label{fig:ex1}
\end{figure}

Figure \ref{fig:ex1} shows the numerical simulation of the Euler-Bernoulli beam with initial conditions $v_0 = 1$, $\xi_0  = 0$ and $z_0 = 0$. It can be seen that the regulation error diminishes very rapidly. In the simulation the spatial derivatives were approximated by finite differences with grid size $0.05$.

\section{Conclusions} \label{sec:conclusions}

We considered robust regulation of impedance passive port-Hamiltonian systems of arbitrary order and showed that a controller satisfying the $\G$-conditions is robust. The robustness result not only holds for impedance passive port-Hamiltonian systems but for any boundary control and observation system satisfying Assumption \ref{ass:standing}. We also presented a simple, minimal order controller structure that satisfies the $\G$-conditions and showed that it stabilizes the closed-loop system, thus solving the robust output regulation problem. The theory was illustrated with an example where we implemented such a controller for a one-dimensional Euler-Bernoulli beam with boundary controls and boundary observations.

\appendix

\section{Admissibility of the observation operator}

\begin{lemma}
Consider a port-Hamiltonian system $(\A, \B, \C)$ as in Definition \ref{def:pHs} and assume that the operator $\B$ is such that $W_B \Sigma W_B^* > 0$. Then the observation operator $\C$ is infinite-time admissible for the semigroup $T_A$ generated by $A = \A|_{\N(\B)}$.
\begin{proof}
Consider the classical solution $x(t) = T_A(t)x_0$ of $\dot{x}(t) = Ax(t)$, $x(0) = x_0 \in \D(A)$ and recall the estimate that was mentioned in the proof of Lemma \ref{lem:wbkwc}:
\begin{equation}
\rp \langle Ax,x \rangle_X \leq \rp\langle f_\partial, e_\partial \rangle_{\mathbb{C}^{nN}}.
\label{eq:Aipest}
\end{equation}
Since $x \in \D(A)$, we have that $\B x = 0$, i.e., $(f_\partial, e_\partial)^T \in \N(W_B)$. As $W_B\Sigma W_B^* > 0$, \cite[Lem. A.1]{LeGZwa05} implies that we may write $W_B = S[I+V_B \quad I-V_B]$ where $S$ is invertible and $V_B$ is square satisfying $V_B^*V_B < I$. Furthermore, as $(f_\partial, e_\partial)^T \in \N(W_B)$, by \cite[Lem. A.2]{LeGZwa05} we may write
\begin{equation}
	\left[\begin{array}{c} f_{\partial}\\ e_{\partial} \end{array}
\right] = \left[\begin{array}{c} I- V_B \\ -I-V_B \end{array} \right]
\ell
\label{eq:kerB}
\end{equation}
for some $\ell \in \mathbb{C}^{nN}$. Let us define the output as $y = \C x$ and write $W_C = [C_1, C_2]$ with $C_{1,2}$ square. We have
$$
  \left[\begin{array}{c} 0 \\y \end{array} \right] =
\left[\begin{array}{c} W_B \\ W_C \end{array}
\right]\left[\begin{array}{c} f_{\partial}\\ e_{\partial} \end{array}
\right]  = \left[\begin{array}{c} 0 \\ C_1(I-V_B) - C_2(I+V_B)
\end{array} \right] \ell
$$
for some $\ell \in \mathbb{C}^{nN}$. Since $\N(W_B)\cap \N(W_C) = \{0\}$, it follows from the above that the square matrix $R:=C_1(I-V_B)-C_2(I+V_B)$ is invertible. Now using the estimate \eqref{eq:Aipest} together with \eqref{eq:kerB} we obtain
$$
\begin{aligned}
\frac{d}{dt}||x(t)||_X^2 & = 2\rp\langle Ax,x \rangle_X \leq 2\rp\langle f_\partial, e_\partial \rangle_{\mathbb{C}^{nN}}.\\
& = \ell^*(-2I + 2V_B^*V_B)\ell \\
& = y^*R^{-*}(-2I + 2V_B^*V_B)R^{-1}y \\
& \leq -m||y||_{\mathbb{C}^{nN}}^2,
\end{aligned}
$$
for some $m > 0$ as $V_B^*V_B < I$. Integrating both sides over $[0, \tau]$ and using $y(t) = \C T_A(t)x_0$ yields
$$
||x(\tau)||_X^2 - ||x_0||_X^2 \leq -m\int_0^\tau ||\C T_A(t)x_0||_{\mathbb{C}^{nN}}^2dt.
$$
Letting $\tau \to \infty$, we have $||x(\tau)||_X^2 \to 0$ as $T_A$ is exponentially stable by part \ref{wbkwc:1}) of Lemma \ref{lem:wbkwc}, and we obtain
$$
	\int_0^\infty ||\C T_A(t)x_0||_{\mathbb{C}^{nN}}^2dt \leq \frac{1}{m	}||x_0||_X^2, 
$$
which concludes the proof.
\end{proof}
\label{lem:admis}
\end{lemma}

\bibliographystyle{plain}

\begin{thebibliography}{10}

\bibitem{Aug_Phd}
B.~Augner.
\newblock {\em Stabilization of Infinite-Dimensional Port-{H}amiltonian Systems
  via Dissipative Boundary Feedback}.
\newblock PhD thesis, Bergische Universit{\"a}t Wuppertal, 2016,
  \url{http://nbn-resolving.de/urn/resolver.pl?urn=urn%3Anbn%3Ade%3Ahbz%3A468-20160719-090307-4}.

\bibitem{AugJac14}
B.~Augner and B.~Jacob.
\newblock Stability and stabilization of infinite-dimensional linear
  port-{H}amiltonian systems.
\newblock {\em Evolution Equations and Control Theory}, 3(2):207--229, 2014.

\bibitem{CurZwa_Book}
R.~Curtain and H.~Zwart.
\newblock {\em An Introduction to Infinite-Dimensional Linear Systems Theory},
  volume~21 of {\em Texts in Applied Mathematics}.
\newblock Springer-Verlag, 1995.

\bibitem{Dav_CDC75}
E.~J. Davison.
\newblock Multivariable tuning regulators: The feedforward and robust control
  of a general servomechanism problem.
\newblock {\em Proc. CDC'75}, 1975.

\bibitem{LeGZwa05}
Y.~Le Gorrec, H.~Zwart, and B.~Maschke.
\newblock Dirac structures and boundary control systems associated with
  skew-symmetric differential operators.
\newblock {\em SIAM J. Control Optim.}, 44(5):1864--1892, 2005.

\bibitem{HamPoh00}
T.~H{\"a}m{\"a}l{\"a}inen and S.~Pohjolainen.
\newblock A finite-dimensional robust controller for systems in the
  {CD}-algebra.
\newblock {\em IEEE Trans. Automat. Control}, 45(3):421--431, 2000.

\bibitem{HamPoh02}
T.~H{\"a}m{\"a}l{\"a}inen and S.~Pohjolainen.
\newblock Robust regulation for exponentially stable boundary control systems
  in {H}ilbert space.
\newblock {\em Proc. MMAR'02, Szczecin, Poland, September 2-5}, 2002.

\bibitem{HamPoh10}
T.~H{\"a}m{\"a}l{\"a}inen and S.~Pohjolainen.
\newblock Robust regulation of distributed parameter systems with
  infinite-dimensional exosystems.
\newblock {\em SIAM J. Control Optim.}, 48(8):4846--4873, 2010.

\bibitem{HamPoh11}
T.~H{\"a}m{\"a}l{\"a}inen and S.~Pohjolainen.
\newblock A self-tuning robust regulator for infinite-dimensional systems.
\newblock {\em IEEE Trans. Automat. Control}, 56(9):2116--2127, 2011.

\bibitem{HumPau_ECC16}
J.-P. Humaloja, L.~Paunonen, and S.~Pohjolainen.
\newblock Robust regulation for first-order port-{H}amiltonian systems.
\newblock {\em Proc. EEC'16, Aalborg, Denmark, June 29th - July 1st}, 2016.

\bibitem{HumPau_MTNS16}
J.-P. Humaloja, L.~Paunonen, and S.~Pohjolainen.
\newblock Robust regulation for port-{H}amiltonian systems of even order.
\newblock {\em Proc. MTNS'16, MN, USA, July 12-15}, 2016.

\bibitem{JacZwa_Book}
B.~Jacob and H.~Zwart.
\newblock {\em Linear Port-{H}amiltonian Systems on Infinite dimensional
  Spaces}, volume 223 of {\em Operator Theory: Advances and Applications}.
\newblock Birkh{\"a}user, 2012.

\bibitem{MacLeG17}
A.~Macchelli, Y.~Le Gorrec, H.~Ramirez, and H.~Zwart.
\newblock On the synthesis of boundary control laws for distributed parameter
  port-{H}amiltonian systems.
\newblock {\em IEEE Trans. Automat. Control}, 62(4):1700--1713, 2017.

\bibitem{Pau16}
L.~Paunonen.
\newblock Controller design for robust output regulation of regular linear
  systems.
\newblock {\em IEEE Trans. Automat. Control}, 61(10):2974--2986, 2016.

\bibitem{PauPoh10}
L.~Paunonen and S.~Pohjolainen.
\newblock Internal model theory for distributed parameter systems.
\newblock {\em SIAM J. Control Optim.}, 48(7):4753--4775, 2010.

\bibitem{PauPoh14a}
L.~Paunonen and S.~Pohjolainen.
\newblock The internal model principle for systems with unbounded control and
  observation.
\newblock {\em SIAM J. Control Optim.}, 52(6):3967--4000, 2014.

\bibitem{Pho91}
V.~Q. Ph{\^o}ng.
\newblock The operator equation {$AX - XB = C$} with unbounded operators {$A$}
  and {$B$} and related abstract {C}auchy problems.
\newblock {\em Math. Z.}, 208:467--588, 1991.

\bibitem{RamLeG14}
H.~Ramirez, Y.~Le Gorrec, A.~Macchelli, and H.~Zwart.
\newblock Exponential stabilization of boundary controlled port-{H}amiltonian
  systems with dynamic feedback.
\newblock {\em IEEE Trans. Automat. Control}, 59(10):2849--2855, 2014.

\bibitem{RebWei03}
R.~Rebarber and G.~Weiss.
\newblock Internal model based tracking and disturbance rejection for stable
  well-posed systems.
\newblock {\em Automatica}, 39:1555--1569, 2003.

\bibitem{TucWei_Book}
M.~Tucsnak and G.~Weiss.
\newblock {\em Observation and Control for Operator Semigroups}.
\newblock Advanced Texts. Birkh{\"a}user Verlag AG, 2009.

\bibitem{VilZwa09}
J.A. Villegas, H.~Zwart, Y.~Le Gorrec, and B.~Maschke.
\newblock Exponential stability of a class of boundary control systems.
\newblock {\em IEEE Trans. Automat. Control}, 54(1):142--147, 2009.

\bibitem{VilZwa_CDC05}
J.A. Villegas, H.~Zwart, Y.~Le Gorrec, B.~Maschke, and A.J. van~der Schaft.
\newblock Stability and stabilization of a class of boundary control systems.
\newblock {\em Proc. CDC'05, Seville, Spain, December 12-15}, 2005.

\end{thebibliography}

\end{document}